\documentclass[final,epsfig,amsfont]{article}
\usepackage{latexsym}
\usepackage[dvips]{graphicx}
\newcommand{\eeq}{\end{equation}}
\newcommand{\beq}{\begin{equation}}

\newtheorem{definition}{Definition}
\newtheorem{conjecture}{Conjecture}

\newtheorem{lemma}{Lemma}
\newtheorem{theorem}{Theorem}

\newcommand{\qed}{\rule{2mm}{3mm}}

\begin{document}
\title{On the Attractor of One-Dimensional Infinite Iterated Function Systems}
\author{
Giorgio Mantica \\
Center for Non-linear and Complex Systems,\\
Department of Science and High Technology, \\ University of Insubria, 22100
Como, Italy \\  and \\ I.N.F.N. sezione di Milano,  CNISM unit\`a di Como. }
\date{}
\maketitle

\begin{abstract}
We study the attractor of Iterated Function Systems composed of infinitely many affine, homogeneous maps. In the special case of {\em second generation IFS}, defined herein, we conjecture that the attractor consists of a finite number of non-overlapping intervals. Numerical techniques are described to test this conjecture, and a partial rigorous result in this direction is proven.
\end{abstract}


{\em Keywords: Iterated Function Systems -- Attractors -- Second Generation IFS} \\

\section{Introduction and statement of results}
\label{sec1}

Iterated function systems (IFS) \cite{papmor,hut,diaco,dem,ba2} are
collections of maps $\phi_i : {\bf R}^n \rightarrow
{\bf R}^n$, $i = 1, \ldots, M$, for which there exists a
set ${\mathcal A}$, called the {\em attractor} of the IFS, that solves
the equation
 \begin{equation}
 \label{attra}
    {\mathcal A}=\bigcup_{i=1,\ldots ,M}\;\phi_i({\mathcal A}) :=  \Phi ({\cal A}).
 \end{equation}
Existence and uniqueness of ${\mathcal A}$ can be easily proven to hold for hyperbolic IFS, {\em i.e.} those for which the maps $\phi_i$ are contractive. In this case, the right-hand side of eq. (\ref{attra}) defines an operator $\Phi$ on the set of compact subsets of ${\bf R}^n$ whose fixed point is ${\cal A}$.
Since $\Phi$ is contractive in the Hausdorff metric,
the set ${\cal A}$ can be also found as the limit of the sequence $\Phi^{n}(K_{0})$,
where $K_{0}$ is any non-empty compact set, i.e.
\begin{equation}
{\cal A} = \lim_{n \to \infty} \Phi^{n}(K_{0}).
\label{eq-itera}
\end{equation}

Attractors of Iterated Function Systems feature a rich variety of topological structures, so that their full characterization is far from being fully understood.
Even in the one-dimensional case, the attractor of an IFS can take on quite different forms. Consider in fact the one dimensional IFS composed of affine maps:
\begin{equation}
\label{mappi}
    \phi_i (s) = \delta_i (s - \beta_i) + \beta_i  \;\;  i = 1, \ldots, M ,
\end{equation}
where $\delta_i$ are real numbers between zero and one, called {\em contraction ratios}, while
$\beta_{i}$ are real constants, that geometrically correspond to the fixed points of the maps. Taking just two maps $\phi_1(x) = \delta x$, $\phi_2(x) = \delta x + 1-\delta$, when $\delta \geq 1/2$ the attractor is the full interval $[0,1]$. To the contrary, when $\delta$ is smaller that one half, the attractor is a Cantor set.
Next, consider the set of three maps:
$\phi_1(x) = x/2$,  $\phi_2(x) = x/4 + 1/4$, $\phi_3(x) = x/4 + 3/4$,
suggested to us by Frank Mendivil, who showed that its attractor is composed of a countable set of disjoint intervals accumulating at one.

In more dimensions, the problem of what compact sets appear as attractors of IFS is even more delicate \cite{duvall,hata,kwie}.
Clearly, many of the technical difficulties in this characterization are typical of many dimensional spaces. Far from wishing to attack this problem, in this paper we focus on one--dimensional systems, albeit of a very special kind: the IFS we consider are composed of uncountably many maps. IFS with infinitely many maps have been studied in \cite{urba,fernau,moran}, in the countable case. Here, to construct an uncountable set of maps we generalize the notion of finite, homogeneous IFS following Elton and Yan \cite{elton} (successively studied and refined in \cite{nalgo1,nalgo2,arxiv}) to define a $(\delta,\sigma)$--{\em homogeneous} affine IFS as follows.
\begin{definition}
\label{gafhom}
Let $\sigma$ be a positive Borel probability measure
on ${\bf R}$ whose support is contained in $[0,1]$,
let $\delta$ be a real number in $[0,1)$ and let $\bar{\delta} := 1 - \delta$.  Let the real number $\beta$ parameterize the IFS maps $\phi_\delta(\beta,\cdot)$ as
$$
 \phi_\delta(\beta,s) := \delta s + \bar{\delta} \beta.
$$
The invariant IFS
measure associated with the { affine} $(\delta,\sigma)$--{homogeneous} IFS is the unique probability measure $\mu$ that satisfies
\begin{equation}
\label{bala2}
    \int   f(s) \; d\mu(s)  =    \int d \sigma(\beta) \int d\mu(s) \;
    f (\phi_\delta(\beta,s)),
\end{equation}
for any continuous function $f$.
\end{definition}
General theory \cite{elton,mendiv} guarantees that the measure $\mu$ defined above is unique. It is also termed a {\em balanced} measure.
The usual, finite homogeneous IFS can be obtained by using a point measure $\sigma_p$ in place of $\sigma$ in eq. (\ref{bala2}):
\begin{equation}
\label{sigma1}
   \sigma_p = \sum_{j=1}^M \pi_j D_{\beta_j},
\end{equation}
where $D_{x}$ is a unit mass, atomic (Dirac) measure at the point $x$ and $\pi_j$, $j=1,\ldots,M$, $\pi_{i} > 0$,
$\sum_{i} \pi_{i} = 1$, are the usual IFS weights.

In words, Definition \ref{gafhom} means that the set of IFS maps is composed of  affine maps, with homogeneous contraction ratio $\delta$, and fixed points $\beta$ distributed according to the measure $\sigma$.  The invariant measure $\mu$ can also be constructed via the usual ``chaos game'', now generalized to an infinity of maps.
Construct a stochastic process in $X=[0,1]$ via the
following rule: given a point $x \in [0,1]$, choose a value of
$\beta$ at random in $[0,1]$, according to the distribution $\sigma(\beta)$ and apply the function
$\phi_{\delta}(\beta,\cdot)$ to map $x$ into $\phi_{\delta}(\beta,x)$. In so doing, the measure $\mu$ can be found, probability one, by the
Cesaro average of atomic measures at the points $x_j$ of a trajectory of the process: $\frac{1}{n} \sum \delta_{x_j} \rightarrow \mu$.

The properties of the measure $\mu$ as a function of $\sigma$, like {\em e.g.} singular versus absolute continuity have been studied in \cite{nalgo2,arxiv}. Approximation and inverse problems where considered in \cite{complex,gio1,nalgo1}, Jacobi matrix construction in \cite{cap,arxiv}. We now focus on a topological, rather than measure theoretical, problem: the structure of the attractor of such IFS. Nonetheless, find convenient to characterize ${\cal A}$ as the support of the measure $\mu$: ${\cal A} : ={\cal S}_\mu$.

This problem is still very general. Therefore, we further restrict our consideration to a specific class of measures $\sigma$ in Definition \ref{gafhom}: those that are themselves the invariant measure of a finite IFS of the kind (\ref{mappi}). This statement needs to be explained in full detail, to avoid confusions:
we start from a finite IFS (that we call a {\em first-generation IFS}) whose invariant measure we label as $\sigma$. We use this symbol because we successively use such $\sigma$ in eq. (\ref{bala2}) to construct a second, homogeneous $(\delta,\sigma)$-IFS with contraction ratio $\delta$ and distribution of fixed points $\sigma$. In so doing, eq. (\ref{bala2}) provides us with the invariant measure $\mu$ we want to study:
\begin{equation}
\label{iter}
\begin{array}{l}
  \{\delta_i,\beta_i,\pi_i\}_{i=1,\ldots,M} \stackrel{1}{\Longrightarrow} \sigma \\
  \{\sigma, \delta \} \stackrel{2}{\Longrightarrow} \mu.
  \end{array}
\end{equation}
The lines in this scheme describe the sequence of the {\em first} and {\em second generation} IFS that are considered in this paper, and the arrows point to the invariant measures that are generated by the respective IFS's.  Notice that an operator $\Phi$ of the kind (\ref{attra}) can be associated to each IFS: we shall use the same letter $\Phi$ for both, labeling them with the index 1 or 2, when necessary. Our aim will then be to find ${\cal A}$, the fixed point of $\Phi_2$. We will call this latter system a {\em second-generation IFS}.
\begin{definition}
A second-generation IFS is a homogeneous IFS, with contraction ratio $0 < \delta < 1$, whose distribution of fixed points $\sigma$ is the invariant measure of a finite maps IFS.
\label{def-second}
\end{definition}

We will mostly assume in this paper that the convex hull of the support of $\sigma$ and $\mu$ is $[0,1]$. In particular, this requires that zero and one be the fixed points of a map of both first and second generation IFS. Also, the first generation IFS may, or may not, be homogeneous, and typically we will consider it as non--overlapping.

The main result of this paper is a conjecture on the nature of the fixed point of $\Phi_2$:
\begin{conjecture}
The attractor of a second-generation affine, homogeneous IFS with $1> \delta > 0$ and disconnected first-generation IFS is composed of a finite number of non-overlapping intervals.
\label{conjint}
\end{conjecture}

To arrive at this conjecture, in Sect. \ref{sec-gen} we first describe two useful lemmas on the support of a generic IFS, and on the action of the operator $\Phi_2$ on intervals. They permit to derive an algorithm for the actual computation of the attractor $\cal A$, in section \ref{sec-num}. This algorithm converges in a finite number of iterations if and only if the attractor  $\cal A$ verifies conjecture \ref{conjint}. We always observe this fact in our numerical experiments. In section \ref{sec-approx} we approximate $\cal A$ from the outside, via the complement of a finite set of open intervals, explicitly computed. We observe numerically that this approximation is sometimes exact, and typically rather satisfactory. We finally conclude in section \ref{suppo} with a partial result in the way of proving conjecture \ref{conjint}: we prove rigorously that for certain second generation IFS the set $\cal A$ contains at least an interval.

\section{General results on the support of the measure $\mu$.}
\label{sec-gen}
In this section we present two results that will be useful in the next construction of the attractor $\cal A$.
We first quote a general result, that holds for any measure $\sigma$, and not only for those considered in the sequel. It shows that the support of $\mu$ is not too far from that of $\sigma$. To simplify formulae it is convenient here to take $[-1,1]$ as the convex hull of $\sigma$ and $\mu$.
\begin{lemma}
Let $S_{\mu}$, $S_{\sigma}$ be the supports of $\mu$ and $\sigma$. Let $S_{\sigma} \subset [-1,1]$.
Then, for any $\delta > 0$, $S_{\sigma} \subset S_{\mu} \subset B_{2 \delta}(S_{\sigma})$.
\label{lem-one}
\end{lemma}
{\em Proof.} See \cite{nalgo2,arxiv}. In the above, $B_{2 \delta}(S_{\sigma})$ is the $2 \delta$--neighborhood of $S_\sigma$.

The second result considers the images of an interval under a finite number of homogeneous IFS maps.

\begin{lemma}
Consider a subset of IFS maps,
   $  \phi_{j} (x) = \delta  (x - \beta_j) + \beta_j,  $
where $\beta$ takes the set of increasingly ordered values $\{ \beta_j , j \in {\cal J} \}$ of finite cardinality.
Let $I = [A,B]$ be an arbitrary interval and let $l$ be its length: $l := B - A$. Suppose that there exist $k$ and $h$ such that
\begin{equation}
   \beta_{j+1} - \beta_j \leq l  {\delta}/{\bar{\delta}}  \;\; \mbox{for all } \; j = k, \ldots, k+h-1.
   \label{condov1}
\end{equation}
Then, the action of the operator $\Phi$ on $I$ contains an interval:
\begin{equation}\label{overlapChain}
 \Phi(I) \supseteq    \bigcup_{j=k}^{h+k} \phi_j(I) = [\phi_k (A), \phi_{k+h}(B)].
\end{equation}
 \label{lem-dista}
 \end{lemma}

{\em Proof.}
Let $I_j = \phi_j (I):= [a_j,b_j]$, for  $j = k,k+1,\ldots,k+h$.
Observe that $a_j = \phi_j(A)$, $b_j = \phi_j(B)$. We obviously suppose that $k,h>0$ and $B>A$.
Clearly, when $b_j \geq a_{j+1}$ we have that
$I_j \cap I_{j+1} \neq \emptyset$. A simple computation reveals that this is equivalent to
$
 \beta_{j+1} - \beta_j \leq l  {\delta}/{\bar{\delta}}.
$
If this holds for all $j = k, \ldots, k+h-1$, then the intervals $I_j$ form an overlapping chain, and eq (\ref{overlapChain}) holds. \qed

Remark that the above lemma requires that the distances between all successive fixed points between $\beta_k$ and $\beta_{k+h}$ must be smaller than the quantity at r.h.s. of eq. (\ref{condov1}), that is a constant. Therefore, only relative positions matter, and not the location of the $\beta$'s.
Furthermore, letting $d_j = (\beta_{j+1} - \beta_j)/l$, we can rewrite condition (\ref{condov1})  as
\begin{equation}\label{djcondition}
   \delta \geq \frac{d_j}{1+d_j} \;\; \mbox{for all } \; j = k, \ldots, k+h-1,
\end{equation}
that shows that for any choice of $k,k+h  \in {\cal J}$ there is a minimal value of $\delta$ for which condition (\ref{condov1}) holds.

\section{Numerical evaluation of the support of the measure $\mu$.}
\label{sec-num}

Suppose now that the distribution of fixed points $\sigma$ is generated by a non-overlapping IFS with a finite number of maps, of the kind (\ref{mappi}): that is, let us consider a {\em second-generation IFS}, eq. (\ref{iter}). We can devise a numerical algorithm to compute the action of the operator $\Phi_2$ on any interval $I$:
\begin{equation}
\Phi_2(I) = \bigcup_{\beta \in {\cal S}_\sigma} \phi_\delta(\beta,I) = \bigcup_{\beta \in {\cal S}_\sigma} ( \delta (I) + \bar{\delta} \beta ).
\label{equai}
\end{equation}
Clearly, since the support of $\sigma$ is uncountable, the above definition is not amenable of numerical treatment. Nonetheless, we can make use of Lemma \ref{lem-dista} above. In doing this,
we find it convenient to construct a countable set of points in
${\cal S}_\sigma$, the {\em band edges}. In fact, under the conditions specified above, the set $\Phi_1^n([0,1])$ is composed of $M^n$ disjoint intervals, that we can call the bands at iteration $n$. For simplicity, label these intervals as $[a^n_j,b^n_j]$. The extrema of these intervals constitute the set of band edges. Let now $l$ be the length of the interval $I$ in eq. (\ref{equai}). Then, when $b^n_j - a^n_j \leq l  {\delta}/{\bar{\delta}}$ holds, Lemma \ref{lem-dista} implies that we can write
\begin{equation}
 \bigcup_{\beta \in {\cal S}_\sigma \cap [a^n_j,b^n_j]} \phi_{\delta}(\beta,I) = \phi_\delta(a^n_j,I) \bigcup  \phi_\delta(b^n_j,I).
\label{equai2}
\end{equation}
That is, out of the uncountable set of maps corresponding to values of $\beta$ in the $i$-th band at iteration $n$ of $\Phi_1$, just two are enough to compute the image of the interval $I$. We can use this observation as the basis of the following algorithm.
\begin{itemize}
 \item[] {{\bf A1. Computing the action of $\Phi_2$ of an interval $I$} \\
 {\bf Input}: the IFS parameters $\{\delta_i,\beta_i\}_{i=1,\ldots,M}$, the  contraction ratio $\delta$, the interval $I$.
 \\ {\bf Output}: the set $\Phi_2 (I)$ as a finite union of $P$ non--overlapping intervals.}
\item[0:] Compute $ \epsilon := l  {\delta}/{\bar{\delta}}$. Initialize the set of band edges with $n=0$, $J=1$, $[a^0_1,b^0_1]=[0,1]$.  Set $L=0$.
\item[1:] For $j=1$ to $J$ and $i=1$ to $M$: Compute the next iteration intervals $\phi_i([a^n_j,b^n_j])$.
\item[2:]  Update $n$ to $n+1$ and $J$ to $MJ$. Set $Z=0$.
\item [3:] For $j=1$ to $J$: Check the inequality $b^n_j-a^n_j \leq \epsilon$. If satisfied, increase $L$ to $L+1$, put $[a^n_j,b^n_j]$ in a list of final points: $[\alpha_L,\beta_L] = [a^n_j,b^n_j]$ and remove it from the list of band edges. Else, increase $Z$ to $Z+1$.
\item [4:] Control. If $Z>0$ set $J=Z$ and loop back to [1]. Else, when $Z=0$ all band edges have been put in the final list, continue.
\item [5:] For $l=1$ to $L$: Compute the interval $\phi(\alpha_l,I) \bigcup \phi(\beta_l,I)$.
\item [6:] By considering intersections, reduce the union of all the intervals in [5] to a sequence of ordered, non intersecting intervals. Compute their cardinality $P$.
\end{itemize}

Observe now that the length $b^n_j-a^n_j$ is certainly less than $\max_i\{\delta_i\}^n$, so that the procedure certainly stops in a finite number of steps.

We now want to apply Algorithm A1 to compute the attractor ${\cal A}$ via eq. (\ref{eq-itera}), starting from the convex hull $K_0 = [0,1]$: $K_n = \Phi_2^{n}([0,1])$. From what demonstrated above, $K^n$ is the union of a finite number of non--overlapping intervals. In the limit, $K_n$ tends (in the Hausdorff metric) to the attractor $\cal A$. It is a matter of experimental observation, that we want to report in this paper, that in all cases we have examined there exists a finite power $\bar{n}$ at which the limit is attained: $\Phi_2^{\bar{n}}([0,1]) = \Phi_2^{\bar{n}-1}([0,1])$. This can be numerically verified by a second algorithm
\begin{itemize}
 \item[] {{\bf A2. Computing the action of $\Phi_2^n$ on $[0,1]$}. \\
 {\bf Input}: the IFS parameters $\{\delta_i,\beta_i\}_{i=1,\ldots,M}$, the contraction ratio $\delta$, the value $n$.
 \\ {\bf Output}: the set $\Phi_2^n ([0,1])$ as a finite union of $Q$ non--overlapping intervals.}
\item[0:] Set $m=0$, $J=1$. Initialize the set $\Phi_2^0 ([0,1])$ to contain the sole interval $[\alpha_1,\beta_1]=[0,1]$.
    \item [1:] Set $L=0$.
    \item [2:] For $j=1$ to $J$: Apply algorithm A1 to compute $\Phi_2([\alpha_j,\beta_j])$, where $[\alpha_j,\beta_j]$ is the $j$-th item in the list $\Phi_2^m([0,1])$.  Add the $P_j$ resulting intervals to a work list of new intervals. Update $L$ to $L+P_j$.
\item [3:] By considering intersections, reduce the union of all the $L$ new intervals computed in [2] to a sequence of ordered, non intersecting intervals, and store it into the list $\Phi_2^{m+1} ([0,1])$. Compute their cardinality $Q$, update $J$ to $Q$.
    \item[4:] Control. If the computed set $\Phi_2^{m+1}([0,1])$ is equal to $\Phi_2^{m}([0,1])$, or if $m+1=n$ stop. Else,  increase $m$ by one and loop back to [1].
\end{itemize}

\begin{figure}[!h]
\begin{center}
\includegraphics[width=.6\textwidth, angle = -90]{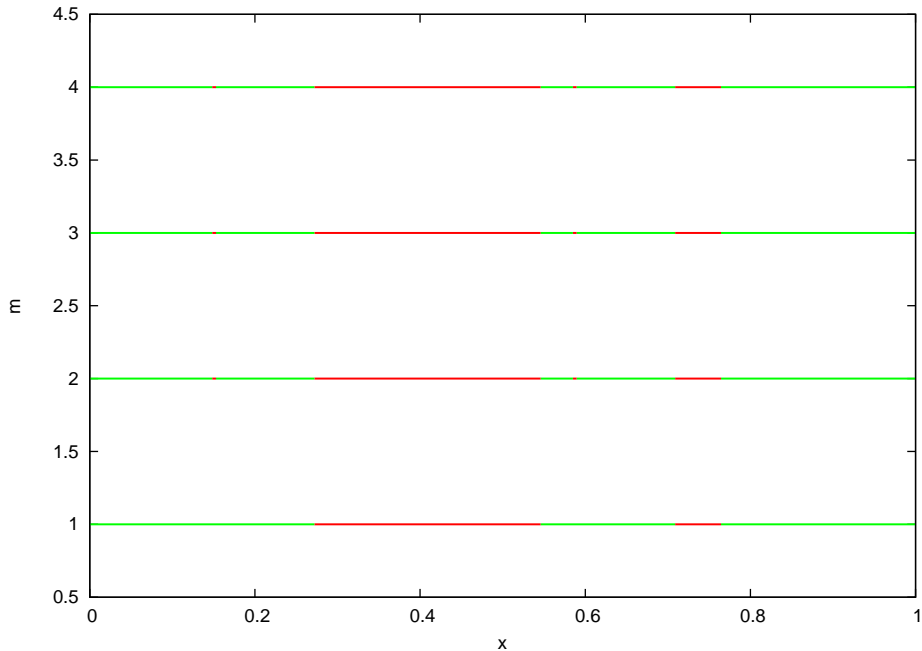}
\caption{Repeated action of the operator $\Phi_2$ on the interval $[0,1]$ for the IFS described in the text. Iteration number is $m$, $\Phi_2^m([0,1])$ is drawn in green and its complement in red. Already at $m=2$ we find ${\cal A} = \Phi_2^m([0,1])$.}
\end{center}
\label{fig-good1}
\end{figure}

As an example of a typical situation, let us now show the application of algorithm A2 to the second-generation IFS given by the two maps $\phi_1(x)=x/5$, $\phi_1(x)=2x/5+3/5$, and $\delta=0.085$. In Figure \ref{fig-good1} we plot the successive iterations $\Phi_2^m([0,1])$ for $m=1$ to $m=4$. We observe that these sets coincide for all $m$ larger than two. Therefore, the support of this measure consists of the union of a finite number of disjoint intervals, in this case five: observe in fact that two tiny gaps also appear, in addition to the two larger ones.

It is remarkable that the same behavior has been found in all the numerical experiments we have carried out. In other words, one might conjecture that {\em the support of a second-generation affine, homogeneous IFS with $\delta > 0$ and disconnected first generation IFS is composed of a finite number of non-overlapping intervals.}

\begin{figure}[!h]
\begin{center}
\includegraphics[width=.7\textwidth, angle = -90]{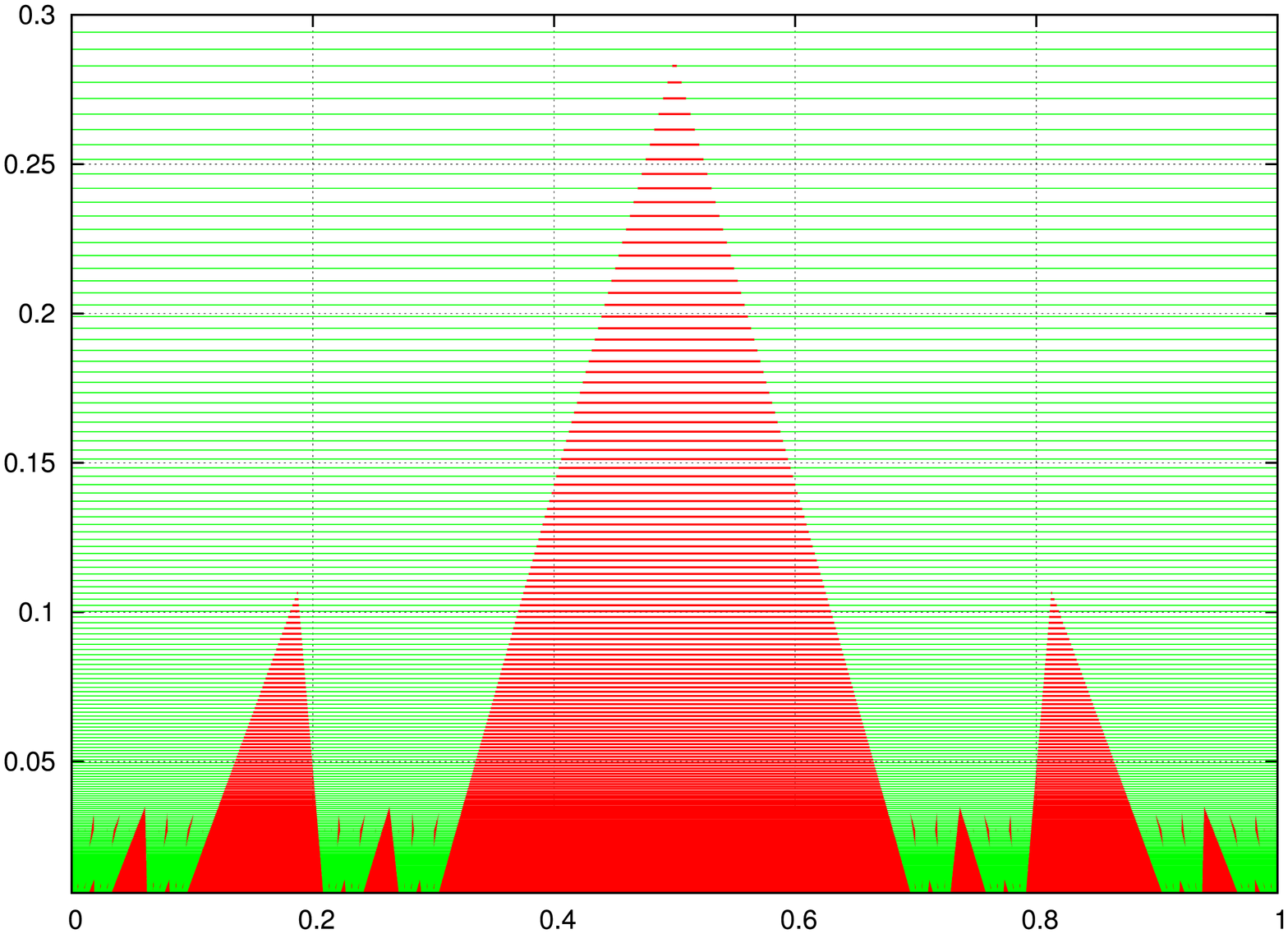}
\caption{Support of the IFS described in the text (green) and gaps (red) stacked one upon the other for varying values of $\delta$ (vertical axis).}
\label{linearScale}
\end{center}
\end{figure}

\begin{figure}[!h]
\begin{center}
\includegraphics[width=.7\textwidth, angle = -90]{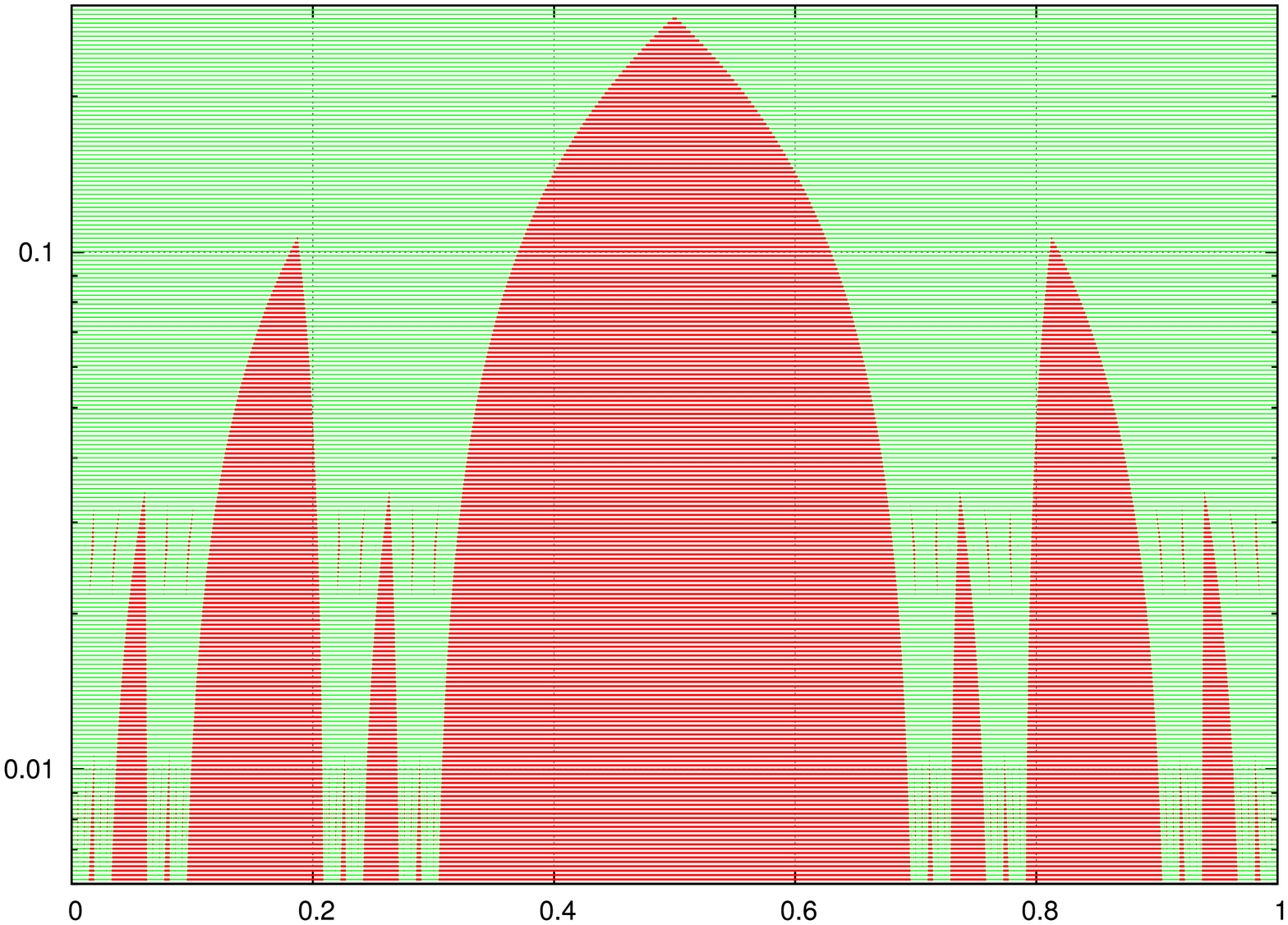}
\caption{Same as Fig. \ref{linearScale}, with vertical axis in logarithmic scale.}
\label{logarithmicScale}
\end{center}
\end{figure}

Figure \ref{linearScale} gives a further illustration of this fact. The basic IFS is generated by the maps $\phi_1(x)=3 x /10$, $\phi_1(x)=3x/10+7/10$, and we let the second-generation contraction ratio $\delta$ vary between $\delta=.006$ and $\delta = .1$.
It is immediately observed that the support is composed of a finite number of intervals, for any finite value of $\delta$, and that this number increases as $\delta$ diminishes. This is perfectly understandable since, as $\delta$ tends to zero, the measure $\mu$ tends to the measure $\sigma$. We will further develop this observation in the next section. The same phenomenon is more evident in Fig. \ref{logarithmicScale}, where the variation of $\delta$ is reported in logarithmic scale.

\section{Approximating the attractor.}
\label{sec-approx}

\begin{figure}[!h]
\begin{center}\label{enneEpsilon}
\includegraphics[width=.7\textwidth, angle = -90]{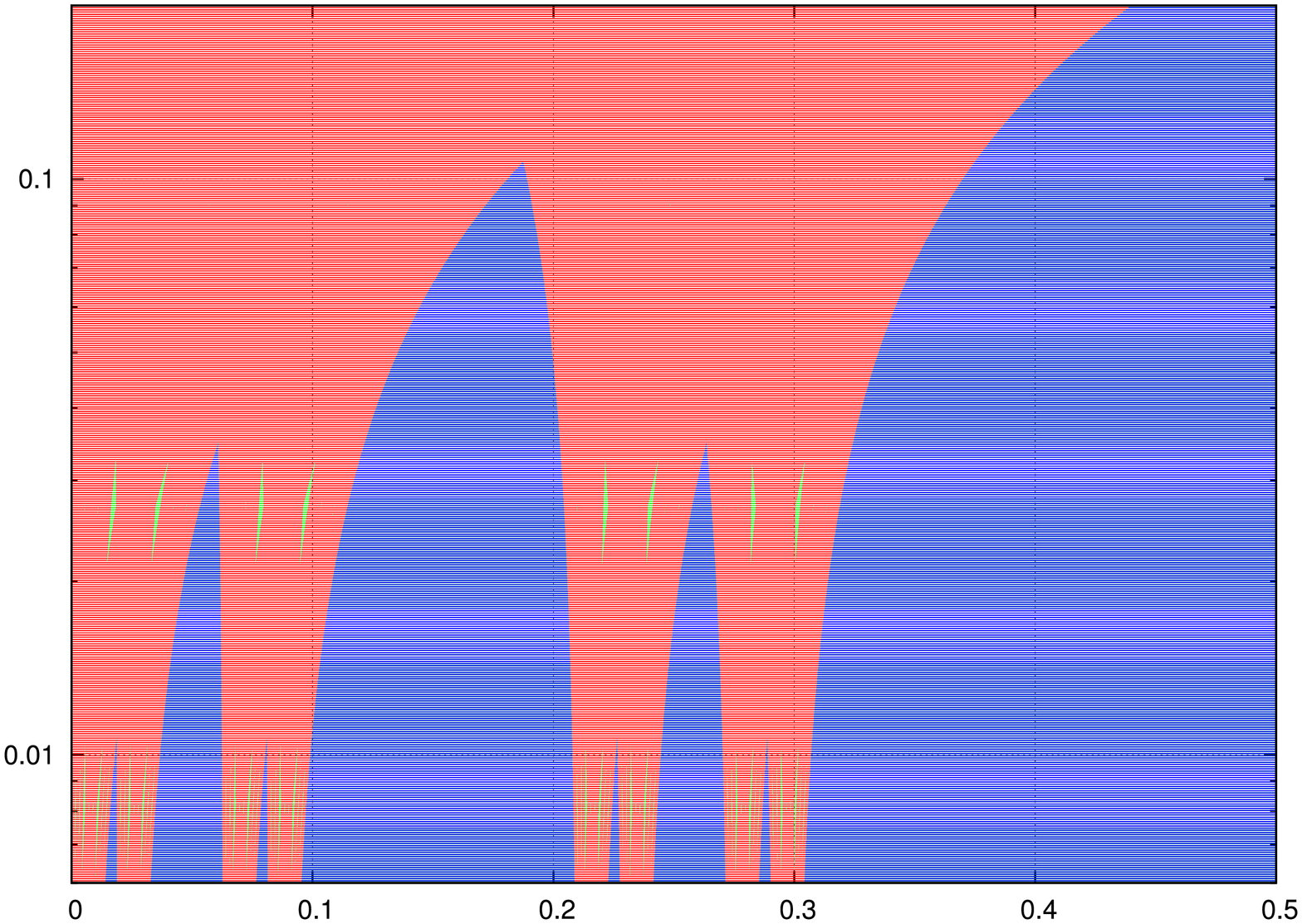}
\caption{Magnification of a segment of figure \ref{logarithmicScale}, with $N_0$ (blue line), $\cal A$ (red) and $G - N_0$ (green) as a function of $\delta$.}
\end{center}
\end{figure}

Observe the detail of the gaps in the support of $\mu$ in Fig. \ref{logarithmicScale}. Most of this structure can be explained by a refinement an analysis similar to that of Lemma \ref{lem-one}. Letting $f$ in eq. (\ref{bala2}) be the characteristic function of $B_\epsilon(x)$, the ball of radius $\epsilon$ centered at $x$, we get the formula
\begin{equation} \label{recall}
\mu(B_{\epsilon} (x)) = \int d \mu(s) \sigma \bigg ( B_{\epsilon/ \overline{\delta}}\bigg ( \frac{x - \delta s}{\overline{\delta}}\bigg )\bigg ).
\end{equation}
Next, let's take into account the fact that the support of $\sigma$ is enclosed in $[0,1]$.
This tell us that whenever $x$ is such that the intersection of $ B_{\epsilon/ \overline{\delta}} \big ( \frac{x - \delta s}{\overline{\delta}}\big)$ with $S_{\sigma}$ is empty for all $s \in [0,1]$, then $\mu(B_{\epsilon} (x)) = 0$. Formally, we can write this condition as:
\begin{equation}\label{firstRewrite}
\bigcup_{s\in [0,1]} \bigg [ \frac{x - \delta s}{\overline{\delta}} - \frac{\epsilon}{\overline{\delta}},
  \frac{x - \delta s}{\overline{\delta}} + \frac{\epsilon}{\overline{\delta}} \bigg ] \cap S_{\sigma} = \emptyset
 \Rightarrow \mu(B_{\epsilon} (x)) = 0.
\end{equation}
The union of intervals at l.h.s. can be easily computed
so that we can rewrite (\ref{firstRewrite}) as
\begin{equation}
\bigg [ \frac{x -\delta}{\overline{\delta}} - \frac{\epsilon}{\overline{\delta}} , \frac{x}{\overline{\delta}} + \frac{\epsilon}{\overline{\delta}}  \bigg] \cap S_\sigma = \emptyset \Rightarrow \mu(B_{\epsilon} (x)) = 0,
\end{equation}
i.e.
\begin{equation}
[x - \delta - \epsilon, x + \epsilon] \cap \overline{\delta} S_{\sigma} = \emptyset \Rightarrow \mu(B_{\epsilon} (x)) = 0.
\label{eq-neps}
\end{equation}
Define now $N_{\epsilon}$ precisely as the set of points that verify the l.h.s. of condition (\ref{eq-neps}). It is easily seen that $N_{\epsilon}$ is the union of a finite number of open intervals, for any $\epsilon$ including zero. From eq. (\ref{eq-neps}) it follows that $N_{\epsilon}$ is enclosed in $G$, the complement of the spectrum:
\begin{equation}
N_{\epsilon} \subseteq G := \overline{S}_\mu.
\label{eq-neps2}
\end{equation}
This latter set, $G$, is the set of {\em gaps}, a finite or countable set of intervals. Therefore, the complementary set of $N_{\epsilon}$ provides an estimate of $S_\mu$ from the outside:
\begin{equation}
 \overline{N}_{\epsilon} \supseteq S_\mu
\label{eq-neps3}
\end{equation}
All the above is true for any $\epsilon$, so that we may let its value tend to zero. It turns out that it is relatively easy to compute numerically $N_{\epsilon}$ and $N_0$ using the same ideas employed in algorithm A1.
Let us therefore examine the nature of the set of gaps in the support of $\mu$, as the union of $N_0$ and a residual set $G - N_0$. We have done this in the same numerical example presented in Figures \ref{linearScale} and \ref{logarithmicScale}.  This is shown in Fig. \ref{enneEpsilon}: $N_0$ is drawn in blue, $\cal A$ in red and $G - N_0$ in green.  We observe that most of $G$ is accounted for by $N_0$, while a non--empty difference $G - N_0$ is observed only when $\delta$ takes values in specific ranges. Therefore, the approximation of eq. (\ref{eq-neps3}) is rather good. It remains therefore to be proven rigorously that $G - N_0$ consists of a finite number of intervals, as experimentally observed. But this cannot be done with the technique of this section. We therefore move on to a deeper approach.

\section{A rigorous result} \label{suppo}

We want to prove now that the support of the measure $\mu$
contains an interval at least.
In this perspective, it is best to consider a Fourier space
representation. Take therefore $f(x)=e^{-i y x}$ in eq. (\ref{bala2}) and
use the notation
    \begin{equation} \label{ft1}
         \hat{\nu}(y) := \int d\nu(x) e^{-i y x}
     \end{equation}
to indicate the Fourier transform of an arbitrary measure $\nu$, to get the well known relation \cite{elton,nalgo2}
  \begin{equation}\label{ft2}
    \widehat{\mu}(y) = \hat{\mu}(\delta y) \cdot     \hat{\sigma}(\bar{\delta} y)
    \end{equation}
that links the Fourier transforms of $\sigma$ and $\mu$.
This implies the following:
\begin{lemma}
 \label{lem-fprod}
The invariant measure  $\mu$ of an affine, homogeneous $(\delta,\sigma)$-IFS is an infinite convolution product of rescaled copies of
the measure $\sigma$,
\begin{equation}
\mu(y) = \sigma(y/\delta^{0} \overline{\delta}) * \sigma(y/\delta^{1} \overline{\delta} y) * \sigma(y/\delta^{2} \overline{\delta} y) * \cdots
 \label{fancazzo1}
\end{equation}
  \end{lemma}
{\em Proof.}
By iterating equation (\ref{ft2}) one obtains the Fourier transform $\mu$ in the form of the infinite product
 \begin{equation}   \widehat{\mu}(y) = \prod_{j=0}^{\infty}     \hat{\sigma}(\delta^j \bar{\delta} y).
     \end{equation}\label{ftr3}
Using the basic property $
\widehat{f*g} = \widehat{f} \cdot \widehat{g}
$  and bijectivity  of the Fourier transform
one obtains the thesis. \qed

The fact that $\mu$ is an infinite convolution product of rescaled copies of
the measure $\sigma$, is known when $\sigma$ is a Bernoulli
measure and $\widehat{\mu}(y)$ is an infinite product of trigonometric
functions \cite{berno}. The above lemma extends this fact to the most general situation.

Observe now that the convolution of two measures $\mu$ and $\nu$ is the measure $\lambda :=\mu * \nu$ such that, for any continuous or measurable function $f$,
\begin{equation}
\int f(x) d \lambda = \int \int f (x + y) d \mu (x) d\mu (y).
\end{equation}
If we choose $f(x) = \chi_{E}(x)$ we get
\begin{equation}
(\mu * \nu) (E)= \int \int \chi_{E} (x + y) d \mu (x) d\nu (y),
\end{equation}
so that by the previous formula
\begin{equation}
S_{\mu*\nu} = { \{ z = x + y, x \in S_{\mu}, y \in S_{\nu} \}}.
\label{fancaz2}
\end{equation}

Let us now consider the case of the invariant measure $\mu$ of a second iteration IFS introduced above. The support of the measure $\sigma$, $S_\sigma$, Cantor set. Formulae (\ref{fancazzo1}) and (\ref{fancaz2}) then imply that the support of $\mu$ is an infinite sum of Cantor sets:
\begin{equation}
S_\mu =  \sum_{j=0}^\infty \bar{\delta} \delta^j S_\sigma.
 \label{fancazzo2}
\end{equation}
Since $\delta < 1$ and since the support of $\sigma$ is bounded, the above series converge.
Cabrelli, Hare and Molter \cite{cabrelli} have considered finite sums of Cantor sets. By using their theory, we can prove:
\begin{theorem}
Let $\sigma$ be the invariant measure of a two--maps, disconnected IFS with contraction ratios smaller than one--third. Let $\mu$ be the invariant measure of a homogeneous $(\delta,\sigma)$-IFS with contraction ratio $\delta$ and distribution of fixed points $\sigma$. Then, for any $\delta$, the support of $\mu$ contains an interval.
\end{theorem}

{\em Proof.}
Theorem 3.2 in \cite{cabrelli} applies (in particular) to finite sums of $n$ Cantor sets $C_j$, each generated by a two--maps IFS, with contraction ratios larger than a positive lower bound $a$ and smaller than one--third. It predicts that, when
\begin{equation}
 (n-1) a^2/(1-a)^3 + a/(1-a) \geq 1
 \label{fancazzo3}
\end{equation}
the sums of these Cantor sets contains an interval.

To apply this theorem to our case, observe that we can take for $C_j$ the set $\bar{\delta} \delta^j S_\sigma$, that is generated by a finite IFS. Let $a$ be the minimum of the contraction ratios of such IFS. Observe that $a$ is then the same for all $j$.
Furthermore, observe that by truncating the infinite summation (\ref{fancazzo2}) to a finite value $n$, the resulting set $S^n_\mu = \sum_{j=0}^n \bar{\delta} \delta^j S_\sigma$ is enclosed in $S_\mu$.
Since by choosing $n$ large enough one can satisfy equation (\ref{fancazzo3}), $S^n_\mu$ contains an interval by \cite{cabrelli}, and so does $S_\mu$. \qed

It is interesting to remark that Cabrelli et al.'s technique also tells us explicitly what is the interval concerned: when applied to our case, this provides the interval $I = [0, \bar{\delta} \delta^n]$.  Observe that the smaller $\delta$, the smaller this interval. Also observe, in the proof of the above theorem, that $n$ does {\em not} depend on $\delta$, but only on the ``first--generation'' IFS.
It finally also follows that all integer powers $\Phi_2^j (I)$ belong to the support of $\mu$:
$
 \Phi_2^j (I) \subseteq S_\mu.
$
As seen before, they consist of a finite union of disjoint intervals (that can reduce to a single interval).

It is likely that suitably generalizing the techniques of \cite{cabrelli} a more general result than the above can be proven. We leave this for further investigation.

{\large \bf Acknowledgements}
Research funded by MIUR-PRIN project {\em Nonlinearity and disorder in classical and quantum transport processes}.

\end{document}